\documentclass[12pt,fleqn]{article}
\usepackage{amsmath}

\usepackage[dvips]{color}
\definecolor{dgreen}{rgb}{0,.8,.3}
\definecolor{lblue}{rgb}{.2,.3,.7}
\input amssym.def
\begin{document}

\renewcommand{\Box}{\rule{2.2mm}{2.2mm}}
\newcommand{\BOX}{\hfill \Box}

\newtheorem{eg}{Example}[section]
\newtheorem{theorem}{Theorem}[section]
\newtheorem{lemma}{Lemma}[section]
\newtheorem{example}{Example}[section]
\newtheorem{remark}{Remark}[section]
\newtheorem{proposition}{Proposition}[section]
\newtheorem{corollary}{Corollary}[section]
\newtheorem{defn}{Definition}[section]
\newenvironment{case}
    {\left\{\def\arraystretch{1.2}\hskip-\arraycolsep \array{l@{\quad}l}}
    {\endarray\hskip-\arraycolsep\right.}

    \title{Vanishing viscosity limit for viscous magnetohydrodynamic equations with a slip boundary condition}
    \author{Xiaoqiang Xie \thanks{School of Science, National University of Defense Technology, Changsha, 410073,
China, E-mail: xqxie@fudan.edu.cn  }\quad Changmin Li
\thanks{School of Management, Fudan University, Shanghai,
200433, China, E-mail: 071025014@fudan.edu.cn  }}
\date{}
   \maketitle

    \abstract{We consider the evolutionary MHD systems,
     and study the the regularity and  vanishing viscosity limit of the
     3-D viscous
     system in a class of bounded domains with a slip boundary condition. We derive the convergence is in $H^{2k+1}$,
    for $k\geq1$, if the initial date holds some sufficient conditions.}
\section{Introduction and results}
Let $\Omega$ be an open bounded domain in $R^3$. We consider the
initial and boundary value problem for the system of viscous MHD
equations
\begin{equation}\left\{\begin{split}&\partial_t\bf{u}-\nu\triangle\bf{u}+(\bf{u}\cdot\nabla)\bf{u}-(\bf{H}\cdot\nabla)\bf{H}+\nabla
p=0\text{ in }\Omega,\\
&\nabla\cdot\bf{u}=0 \text{ in }\Omega,\\
&\partial_t\bf{H}-\mu\triangle\bf{H}+(\bf{u}\cdot\nabla)
\bf{H}-(\bf{H}\cdot\nabla)
\bf{u}=0\text{ in }\Omega,\\
&\nabla\cdot\bf{H}=0 \text{ in }\Omega,\\
&\bf{u}=\bf{u}_0, \bf{H}=\bf{H}_0,\text{ at }t=0,
\end{split}
\right.\label{MS}
\end{equation}
with the following slip without friction boundary conditions
\begin{equation}\bf{u}\cdot\bf{n}=0,
\nabla\times\bf{u}\cdot\bf{\tau}=0,\bf{H}\cdot\bf{n}=0,
\nabla\times\bf{H}\cdot\bf{\tau}=0\text{ on
}\partial\Omega\label{bom}
\end{equation}
where $\nabla\cdot$ and $\nabla\times$ denote the div and curl
operators,  $\bf{n}$ the outward normal vector and $\tau$ any unit
tangential vector of $\partial \Omega$.

The corresponding ideal MHD system is usually equipped with the slip
boundary condition, namely
\begin{equation}\left\{\begin{split}&\partial_t\bf{u}^0+(\bf{u}^0\cdot\nabla)\bf{u}^0-
(\bf{H}^0\cdot\nabla)\bf{H}^0+\nabla
p^0=0\text{ in }\Omega,\\
&\nabla\cdot\bf{u}^0=0 \text{ in }\Omega,\\
&\partial_t\bf{H}^0+(\bf{u}^0\cdot\nabla)\bf{H}^0-(\bf{H}^0\cdot\nabla)\bf{u}^0=0\text{ in }\Omega,\\
&\nabla\cdot\bf{u}^0=0 \text{ in }\Omega,\\
&\bf{u}^0=\bf{u}_0, \bf{H}^0=\bf{H}_0,\text{ at }t=0,
\end{split}
\right.\label{IM}
\end{equation}
\begin{equation}\bf{u}^0\cdot\bf{n}=0, \bf{H}^0\cdot\bf{n}=0 \text{ on }\Omega.\label{boi}
\end{equation}

Our aim is to investigate strong convergence, up to the boundary, of
the solution $(\bf{u},\bf{H})$ of the MHD (\ref{MS}) to the solution
$(\bf{u}^0,\bf{H}^0)$ of the ideal MHD system (\ref{IM}), as
$(\nu,\mu)\rightarrow 0$.

The boundary conditions (\ref{bom}) are a special Navier-type slip
boundary conditions, which allow the fluid to slip at a slip
velocity proportional to the shear stress introduced by Navier
\cite{Nav}. This type of boundary conditions has been used in many
fluid problems (see e.g. \cite{Clo}, \cite{Lia}, \cite{Muc},
\cite{Wat}).

The viscous MHD system in the whole space or with non-slip boundary
conditions has been studied  extensively (see e.g. \cite{Bis},
\cite{Duv}, \cite{Xin}, \cite{Ser}, \cite{Jwu}). The solvability,
regularity of the 3-D viscous MHD system with a slip boundary
condition, we refer to \cite{Yue}.

The issue of vanishing viscosity limits of the Navier-Stokes
equations  is classical and fundamental importance in fluid dynamics
and turbulence theory (see e.g. \cite{Kato2} \cite{Sam}, \cite{Tem},
\cite{XQ1}, \cite{Bei}, \cite{Yue2}, \cite{Con}).

In flat boundary case, the 3-D inviscid limit for solution
$(\bf{u},\bf{H})$ to the slip boundary problem (\ref{MS}) and
(\ref{bom}) has been considered in \cite{Yue}. In \cite{Yue}, they
state the following result. Assume $\nabla\cdot\bf{u}_0=0$,
$\nabla\cdot\bf{H}_0=0$, and  $(\bf{u}_0, \bf{H}_0)\in H^3$ satisfy
the boundary conditions (\ref{bom}). Then, as $(\nu,\mu)\rightarrow
0$,
\begin{equation}(\bf{u}, \bf{H})\rightarrow
(\bf{u}^0, \bf{H}^0)\text{ in }L^p(0,T;H^3(\Omega))\cap
C([0,T];H^2(\Omega)),
\end{equation}
for some $T>0$ and any $p\in [1,+\infty)$, where $(\bf{u}^0,
\bf{H}^0)$ is the solution to the ideal MHD equations (\ref{IM}) and
(\ref{boi}).

It should be noted that the approach encounters great difficulties
for general domains as pointed out by \cite{Bei}. Thus, following
\cite{Bei}, we restrict the problem to a cubic domain
$Q=[0,1]^2_{per}\times(0,1)$ with the boundary conditions on two
opposite faces $z=0$ and $z=1$, and others be assumed periodic,
which was called flat boundary case.

Our approach here is motivated by the idea introduced in \cite{XQ2}
to study the same problems for the Navier-Stokes equations. We prove
the following result.

\begin{theorem}Let the initial data $\bf{u}_0\in V^{2k-1}\cap H^{2k+1}$, $\bf{H}_0\in V^{2k-1}\cap H^{2k+1}$,
 $k\geq 1$. Then there
exist strong solution of the MHD equation (\ref{MS}) and (\ref{bom})
in the "cubic domain" (flat boundary case) on some time interval
$[0,T]$, s.t.
\begin{equation}\begin{split}&\|\bf{u}\|_{L^\infty(0,T;H^{2k+1})}+\|\bf{H}\|_{L^\infty(0,T;H^{2k+1})}\leq
C,\\&
\|\partial_t\bf{u}\|_{L^2(0,T;H^{2k})}+\|\partial_t\bf{H}\|_{L^2(0,T;H^{2k})}\leq
C.\end{split}
\end{equation}And
\begin{equation}(\bf{u}, \bf{H})\rightarrow (\bf{u}^0, \bf{H}^0)
 \text{ in }C([0,T];H^{2k}), \text{ as }(\nu,\mu)\rightarrow
0,
\end{equation}
where $(\bf{u}^0, \bf{H}^0)$ is the unique solution of the ideal MHD
equations (\ref{IM}) and (\ref{boi}).

Further, denoting $\omega^0=\nabla \bf{u}^0$ and $\zeta^0=\nabla
\bf{H}^0$, if
$\|\partial^{2k}_n\omega^0_\tau\|_{L^\infty(0,T;C^2(\partial\Omega))}\leq
C$,
$\|\partial^{2k}_n\zeta^0_\tau\|_{L^\infty(0,T;C^2(\partial\Omega))}\leq
C$,
$\|\partial_t\partial^{2k}_n\omega^0_\tau\|_{L^2(0,T;C^1(\partial\Omega))}\leq
C$,
$\|\partial_t\partial^{2k}_n\zeta^0_\tau\|_{L^2(0,T;C^1(\partial\Omega))}\leq
C$, then
\begin{equation}(\bf{u}, \bf{H})\rightarrow (\bf{u}^0, \bf{H}^0)
 \text{ in }C([0,T];H^{2k+1}), \text{ as }(\nu,\mu)\rightarrow
0.
\end{equation}
\end{theorem}

The paper is organized as follows. Some tools are drawn in section
\ref{section2}. A priori estimates to the MHD systems are given in
section \ref{section3}. The results of vanishing viscosity limit and
the convergence rate are presented in section \ref{section4}.

\section{Notations and preliminaries}\label{section2}

Throughout the rest of this paper, denote by
$\bf{v}_\tau=\bf{v}\cdot\tau$ and $\bf{v}_n=\bf{v}\cdot\bf{n}$ on
the boundary $\partial\Omega$. For the flat boundary case,
$\bf{v}\cdot\bf{n}=0$ and $\nabla\times\bf{v}=0$ are equivalent to
$\bf{v}_n=0$ and $\partial_n\bf{v}_\tau=0$ on $\partial Q$. And
$\partial Q=\{(x,y,z);z=0,z=1\}\cap \overline{Q}$. For convenience,
$\Omega$ and $Q$ may be omitted when we write these spaces without
confusion.

We begin our analysis with a formula of integration by parts.

Let $\Omega$ be a regular open, bounded set in $R^3$. Then, for
sufficiently regular vector fields $\bf{v}$,
\begin{equation}-\int_\Omega\triangle\bf{v}\cdot\bf{v}dx
=\|\nabla\bf{v}\|_{L^2}^2-\int_{\partial\Omega}\partial_{\bf{n}}\bf{v}\cdot\bf{v}d\sigma.
\end{equation}

It is easily  shown that if $\bf{v}$ is, sufficiently regular,
vector fields in a flat boundary domain then
\begin{equation}\partial_{\bf{n}}\bf{v}\cdot\bf{v}=\partial_n\bf{v}_\tau\cdot\bf{v}_\tau+
\partial_n\bf{v}_n\cdot\bf{v}_n.
\end{equation}

It follows that  $\partial_{\bf{n}}\bf{v}\cdot\bf{v}$ vanishes on
the boundary if either of the following conditions is satisfied,
\begin{equation}\begin{split}&(a)\quad\bf{v}\cdot\bf{n}=0, \nabla\times\bf{v}\times\bf{n}=0\text{ on
}\partial\Omega,\\&(b)\quad \bf{v}\times\bf{n}=0\text{ on
}\partial\Omega, \nabla\cdot\bf{u}=0\text{ in
}\Omega.\label{e0}\end{split}
\end{equation}

To study  functions with either of above boundary conditions, we
introduce series of function sets.

Let
$$H=\{\bf{v}\in H^1;\nabla\cdot\bf{v}=0\text{ in }\Omega\},$$
$$V^{-1}=\{\bf{v}\in H; \bf{v}_n=0\text{ on }\partial\Omega\},$$
$$V^{0}=\{\bf{v}\in H; \bf{v}_\tau=0\text{ on }\partial\Omega\},
$$
$$V^{2k}=\{\bf{v}\in H^{2k+1};\partial_n^{2j}\bf{v}\in V^0, j=0,1,\cdots,k\},$$
$$V^{2k+1}=\{\bf{v}\in H^{2k+2};\bf{v}\in V^{-1}, \partial_n^{2j+1}\bf{v}\in V^0, j=0,1,\cdots,k\}.$$

Then, the following propositions are easily obtained
\begin{proposition}\label{pro0}
Let $k\geq 1$. Then $V^k\subseteq V^{k-2}$.
\end{proposition}
\begin{proposition}\label{pro1}Let $k\geq -1$,
and $\bf{v}\in V^{k}$. Then $\partial_n^{j}\bf{v}\in V^{k-j}$,
$j=0,1,\cdots,k+1.$\end{proposition}

 Rewrite (\ref{e0}) with the new
notations,
\begin{lemma}\label{lem0}Let $k\geq 0$ and $\bf{v}\in V^k$. Then
$\partial_n\bf{v}\cdot\bf{v}=0$ on $\partial\Omega.$
\end{lemma}

It should be considered that when $\bf{v}$ is not in $V^0$ or
$V^{1}$. For energy estimates, we construct a boundary layer to fill
the gap.
\begin{lemma}In the flat boundary case, assume $\|\bf{h}_\tau\|_{C^{1}(\partial\Omega)}\leq
C$ for $k\geq 1$. Then, for any $\varepsilon<<1$, there is a
$\bf{v}^\varepsilon\in V^{2k-1}$,
$\chi^\varepsilon=\nabla\times\bf{v}^\varepsilon$,  such that
$\nabla\times\bf{v}^\varepsilon\equiv0$ as
$\varepsilon^{\frac{1}{2}}\leq z\leq 1-\varepsilon^{\frac{1}{2}}$,
furthermore,
\begin{equation}\begin{split}&
\bf{\chi}_\tau^\varepsilon\in C^{2k+1}(\overline{\Omega}),
\bf{\chi}^\varepsilon_n\in C^{2k}(\overline{\Omega}),\\&
\partial_n^{2k}\bf{\chi}^\varepsilon_\tau=\bf{h}_\tau, \partial_n^{2k}\bf{\chi}^\varepsilon_n=0\text{ on
}\partial\Omega, \\&
\|z^i(1-z)^i\partial^{2k+1}_n\bf{\chi}^\varepsilon_\tau\|_{L^{p}}\leq
C\varepsilon^{\frac{1}{2p}+\frac{i-1}{2}},\\
&\|z^i(1-z)^i\partial^{2k}_n\bf{\chi}^\varepsilon_n\|_{L^{p}}\leq
C\varepsilon^{\frac{1}{2p}+\frac{i+1}{2}},\\
&\|\partial^{2k}_n\bf{\chi}^\varepsilon\|_{L^{p}}\leq
C\varepsilon^{\frac{1}{2p}+\frac{i}{2}},\end{split}
\end{equation}for $i\in R^+$, $1\leq p<+\infty$.\label{lem3}
\end{lemma}
\textbf{Proof.} It's trivial to find a function $\varphi(z)\in
C^{1}[0,\infty)$, s.t.
\begin{equation}\left\{\begin{split}
&\varphi(z)=1 \text{ at }z=0,\\
&\varphi(z)=0\text{ at }z\geq 1,\\
&\int_0^1F^{j}(\varphi)(s)ds=0, j=0,1,\cdots,2k-1,\end{split}\right.
\end{equation}
where $F$ is an integrate operator from $C[0,\infty)$ to
$C^1[0,\infty)$, and $F(f)(z)=\int_0^zf(s)ds$, $F^0(f)=f$,
$F^j=F(F^{j-1})$, $j\geq 1$.

Denote by
$\varphi^\varepsilon(z)=\varphi(\frac{z}{\varepsilon^{\frac{1}{2}}})$.
Then,
$$\|z^i\partial_z^j\varphi^\varepsilon\|_{L^p}\leq
C\varepsilon^{\frac{1}{2p}+\frac{i-j}{2}}\text{ for }i\in R^+, j\leq
2, 1\leq p\leq +\infty.$$

Set
$\psi^\varepsilon(z)_\tau=\bf{h}_\tau(0)\varphi^\varepsilon(z)+\bf{h}_\tau(1)\varphi^\varepsilon(1-z)$,
and
$\psi^\varepsilon_3=-\int_0^z\nabla_\tau\cdot\psi_\tau(x,y,s)ds$. It
follows that
$$ \nabla\cdot\psi^\varepsilon=0\text{ in }\Omega,
\psi^\varepsilon_\tau=\bf{h}_\tau \text{ on }\partial\Omega.$$

Next, set $$\chi^\varepsilon=F^{2k}(\psi^\varepsilon).$$ Since
$F^{j}(\varphi^\varepsilon)=0$ on $\partial\Omega$, for
$j=1,2,\cdots,2k$, it follows that
$\partial^{j}_z\chi^\varepsilon_\tau=0$, $\partial^{j+1}_z\chi_3=0$
on $\partial \Omega$, for $j=0,1,\cdots,2k-1$. Furthermore,
$\nabla\cdot\chi^\varepsilon=0$ in $\Omega$. In  other words,
$\chi^\varepsilon\in V^{2k-2}$. Therefore, $\int\chi_3=0$.

Finally, let $\zeta^\varepsilon$ satisfy the following equations
\begin{equation}\left\{\begin{split}&
-\triangle\zeta^\varepsilon=\chi^\varepsilon\text{ in }\Omega,\\
&\zeta^\varepsilon_\tau=0, \partial_z\zeta^\varepsilon_3=0\text{ on
}\partial\Omega.
\end{split}\right.\label{e1}
\end{equation}
The necessary condition  $\int\chi_3=0$ of existence  holds by
classical elliptic theories. Applying div to equation (\ref{e1}),
together with div$\zeta^\varepsilon=0$ on $\partial\Omega$, then
$\nabla\cdot\zeta^\varepsilon=0$ in $\Omega$.

Set $\bf{v}^\varepsilon=\nabla\times\zeta^\varepsilon$ and notice
that
$\nabla\times\nabla\zeta^\varepsilon=-\triangle\chi^\varepsilon$,
then the proof is completed after a simple calculation. \Box

Now, we derive some results of  nonlinearities.
\begin{lemma}$(\bf{u}\cdot\nabla)\bf{v}$ is normal to the boundary,
if either of the following conditions holds
\begin{equation}\begin{split}&(a)\qquad\bf{u}\in V^0, \bf{v}\in V^1,\\
&(b)\qquad\bf{u}\in V^{-1}, \bf{v}\in V^0.
\end{split}
\end{equation}\label{lem1}
\end{lemma}

\begin{lemma}Let $j\geq0$, $\bf{u}, \bf{v}\in V^j$.
Then, $(\bf{u}\cdot\nabla)\bf{v}\cdot\bf{n}=0$ on
$\partial\Omega$.\label{lem4}\end{lemma}

The proof is left to the reader.

\begin{theorem}Let $\bf{u}, \bf{v}\in V^{2k+1}$, $k\geq 1$.
Then, for $0\leq j\leq 2k+1$,
\begin{equation}\left\{\begin{split}
&\partial_n^{j}(\bf{u}\cdot\nabla)\bf{v}\cdot\bf{n}=0 ~\text{ on
}\partial\Omega,~\text{ if
}j\text{ is even},\\
&\partial_n^{j}(\bf{u}\cdot\nabla)\bf{v}\times\bf{n}=0 \text{ on
}\partial\Omega, ~\text{ if }j\text{ is odd}.
\end{split}\right.\end{equation}\label{th0}
\end{theorem}
\textbf{Proof.} It's easily derived,
$$\partial^j_n(\bf{u}\cdot\nabla)\bf{v}=\sum_{i=0}^j
\partial^i_n(\bf{u}\cdot\nabla)\partial^{j-i}_n\bf{v}.$$

If $j$ is odd. $i$ and $j-i$, or $j-i$ and $i$, are
 odd and even, respectively. Recalling
Proposition \ref{pro0}, $\partial^i_n\bf{u}\in V^1$ and
$\partial^{j-i}_n\bf{v}\in V^0$, or $\partial^i_n\bf{u}\in V^0$ and
$\partial^{j-i}_n\bf{v}\in V^1$. It follows Lemma \ref{lem1},
$\partial^i_n(\bf{u}\cdot\nabla)\partial^{j-i}_n\bf{v}\times\bf{n}=0$
on $\partial\Omega$, and the desired result is obtained.

If $j$ is even.  $i$ and $j-i$ are all odd or even. Recalling
proposition \ref{pro0} and Lemma \ref{lem4},
$\partial^i_n(\bf{u}\cdot\nabla)\partial^{j-i}_n\bf{v}\cdot\bf{n}=0$.
And the desired result is obtained. \Box

Finally, denote
$\nabla\times(\bf{u}\cdot\nabla)\bf{v}-(\bf{u}\cdot\nabla)(\nabla\times\bf{v})$
by $F(Du,Dv)$. By appealing to  Theorem \ref{th0}, the following
results can be obtained.

\begin{corollary}
Let $\bf{u}, \bf{v}\in V^{2k+1}$, $k\geq 1$. Then, for $0\leq j\leq
k$,

$$\partial_n^{2j}(\bf{u}\cdot\nabla)(\nabla\times\bf{v})\times\bf{n}=0,
 \partial_n^{2j}F(Du,Dv)\times\bf{n}=0\text{ on
 }\partial\Omega.$$\label{cor}
\end{corollary}

\section{A priori estimates}\label{section3}

Now, we derive formal energy estimates assuming that $\bf{u}_0$,
$\bf{u}$, $\bf{H}_0$, $\bf{H}$ are sufficiently regular. As pointed
out in  \cite{Yue} and \cite{Bei}, the key in studying the vanishing
viscosity limit is to control the vorticity created on the boundary.

Set
$$\omega=\nabla\times \bf{u},
\zeta=\nabla\times \bf{H}.$$

Recalling the boundary conditions (\ref{bom}) together with the
notations introduced in section \ref{section2},
 \begin{equation}\bf{u}\in V^1,\quad \bf{H}\in V^1,\quad
\omega\in V^0,\quad\zeta\in V^0.\label{b0}\end{equation}

By applying the operator curl to both sides of the equation
(\ref{MS}) one gets,
\begin{equation}\begin{split}&\partial_t\omega-\nu\triangle\omega
+(\bf{u}\cdot\nabla)\omega-(\bf{H}\cdot\nabla)\zeta+ F(Du,Du)-
F(DH,DH)=0\text{ in }\Omega.\end{split}\label{om1}
\end{equation}
\begin{equation}\begin{split}&\partial_t\zeta-\mu\triangle\zeta
+(\bf{u}\cdot\nabla)\zeta-(\bf{H}\cdot\nabla)\omega+
F(Du,DH)-F(DH,Du) =0\text{ in }\Omega.\end{split}\label{om2}
\end{equation}

By appealing to Corollary \ref{cor} and (\ref{b0}), one obtains
$\partial_t\omega, \partial_t\zeta, \triangle_\tau\omega,
\triangle_\tau\zeta, \\(\bf{u}\cdot\nabla)\omega,
(\bf{H}\cdot\nabla)\zeta,(\bf{u}\cdot\nabla)\zeta,(\bf{H}\cdot\nabla)\omega,
F(Du,Du), F(DH,DH), F(Du,DH)\text{ and } \\F(DH,Du)$ are all normal
to boundary. Then, by equations (\ref{om1}) and (\ref{om2}),
$\partial^2_n\omega$ and $\partial^2_n\zeta$ are normal to
$\partial\Omega$. It follows that $\omega\in V^2$, $\zeta\in V^2$,
$\bf{u}\in V^3$ and $\bf{H}\in V^3$.

Similarly, applying operator $\partial^{2}_n$ to both sides of
equations (\ref{om1}) and (\ref{om2}). Step by step, the following
result is obtained.
\begin{lemma}Let $\bf{u}$ and $ \bf{H}$  be sufficient regularity. Then,
for $k\in N$,
\begin{equation}\bf{u}\in V^{2k+1}, \bf{H}\in
V^{2k+1}, \omega\in V^{2k}, \zeta\in
V^{2k}.\end{equation}\label{lem2}
\end{lemma}

Applying the operator $\partial^{\alpha}_{x,y,z}$, ($\alpha$ is a
multi-index and $ |\alpha|\leq 2k$) to both sides of equations
(\ref{om1}) and (\ref{om2}), one gets
\[\begin{split}&\partial_t\partial^{\alpha}_{x,y,z}\omega-\nu
\triangle\partial^{\alpha}_{x,y,z}\omega
+(\bf{u}\cdot\nabla)\partial^{\alpha}_{x,y,z}\omega-(\bf{H}\cdot\nabla)\partial^{\alpha}_{x,y,z}\zeta\\&+
\sum_{|\beta|=1,
\beta+\gamma=\alpha}(\partial^{\beta}_{x,y,z}\bf{u}\cdot\nabla)\partial^{\gamma}_{x,y,z}\omega-
\sum_{|\beta|=1,\beta+\gamma=\alpha}(\partial^{\beta}_{x,y,z}\bf{H}\cdot\nabla)\partial^{\gamma}_{x,y,z}\zeta\\&
+\sum_{|\beta|\geq1,
\beta+\gamma=\alpha}(\partial^{\beta}_{x,y,z}\bf{u}\cdot\nabla)\partial^{\gamma}_{x,y,z}\omega-
\sum_{|\beta|\geq1,\beta+\gamma=\alpha}(\partial^{\beta}_{x,y,z}\bf{H}\cdot\nabla)\partial^{\gamma}_{x,y,z}\zeta\\&+
F(D\partial^{\alpha}_{x,y,z}u,Du)+F(Du,D\partial^{\alpha}_{x,y,z}u)-
F(D\partial^{\alpha}_{x,y,z}H,DH)-F(DH,D\partial^{\alpha}_{x,y,z}H)\\
&+\sum_{1\leq|\beta|\leq
|\alpha|-1,\beta+\gamma=\alpha}F(D\partial^{\beta}_{x,y,z}u,D\partial^{\beta}_{x,y,z}u)
-\sum_{1\leq|\beta|\leq
|\alpha|-1,\beta+\gamma=\alpha}F(D\partial^{\beta}_{x,y,z}u,D\partial^{\beta}_{x,y,z}u)
=0,\end{split}
\]
\[\begin{split}&\partial_t\partial^{\alpha}_{x,y,z}\zeta-\mu
\triangle\partial^{\alpha}_{x,y,z}\zeta
+(\bf{u}\cdot\nabla)\partial^{\alpha}_{x,y,z}\zeta-(\bf{H}\cdot\nabla)\partial^{\alpha}_{x,y,z}\omega\\&+
\sum_{|\beta|=1,
\beta+\gamma=\alpha}(\partial^{\beta}_{x,y,z}\bf{u}\cdot\nabla)\partial^{\gamma}_{x,y,z}\zeta-
\sum_{|\beta|=1,\beta+\gamma=\alpha}(\partial^{\beta}_{x,y,z}\bf{H}\cdot\nabla)\partial^{\gamma}_{x,y,z}\omega\\&
+\sum_{|\beta|\geq1,
\beta+\gamma=\alpha}(\partial^{\beta}_{x,y,z}\bf{u}\cdot\nabla)\partial^{\gamma}_{x,y,z}\zeta-
\sum_{|\beta|\geq1,\beta+\gamma=\alpha}(\partial^{\beta}_{x,y,z}\bf{H}\cdot\nabla)\partial^{\gamma}_{x,y,z}\omega\\&+
F(D\partial^{\alpha}_{x,y,z}u,DH)+F(Du,D\partial^{\alpha}_{x,y,z}H)-
F(D\partial^{\alpha}_{x,y,z}H,Du)-F(DH,D\partial^{\alpha}_{x,y,z}u)\\
&+\sum_{1\leq|\beta|\leq
|\alpha|-1,\beta+\gamma=\alpha}F(D\partial^{\beta}_{x,y,z}u,D\partial^{\beta}_{x,y,z}H)
-\sum_{1\leq|\beta|\leq
|\alpha|-1,\beta+\gamma=\alpha}F(D\partial^{\beta}_{x,y,z}H,D\partial^{\beta}_{x,y,z}u)
=0.\end{split}
\]

Next, multiplying both sides of the above equations by
$\partial^\alpha_{x,y,z}\omega$ and $\partial^\alpha_{x,y,z}\zeta$,
respectively,  integrating in $\Omega$, and summing them up. Note
that $\partial^\alpha_{x,y,z}\omega^\nu\in V^0$ and
$\partial^\alpha_{x,y,z}\zeta^\nu\in V^0$ or
$\partial^\alpha_{x,y,z}\omega^\nu\in V^1$ and
$\partial^\alpha_{x,y,z}\zeta^\nu\in V^1$ for $|\alpha|\leq 2k$,
$$\int(\bf{H}\cdot\nabla)\partial_{x,y,z}^\alpha\zeta\cdot\partial_{x,y,z}^\alpha\omega
+\int(\bf{H}\cdot\nabla)\partial_{x,y,z}^\alpha\omega\cdot\partial_{x,y,z}^\alpha\zeta=0,$$
$$\|\nabla\bf{u}\|_{L^\infty}+\|\nabla\bf{H}\|_{L^\infty}\leq
C\|\bf{u}\|_{H^{3}}+C\|\bf{H}\|_{H^3}\leq
C\|\omega\|_{H^{2k}}+C\|\zeta\|_{H^{2k}},$$ and
$$\|\bf{u}\|_{W^{j,4}}+\|\bf{H}\|_{W^{j,4}}
\leq C\|\omega\|_{W^{j-1,4}}+ C\|\zeta\|_{W^{j-1,4}}\leq
C\|\omega\|_{H^{j}}+C\|\zeta\|_{H^{j}}, 2\leq j\leq 2k.$$

  By Lemma \ref{lem0} and summing  up for all $|\alpha|\leq 2k$,  one obtains
$$\frac{1}{2}\frac{d}{dt}(\|\omega\|_{H^{2k}}^2
+\|\zeta\|_{H^{2k}}^2)
+\nu\|\nabla\omega\|_{H^{2k}}^2+\nu\|\nabla\zeta\|_{H^{2k}}^2\leq
C\|\omega\|_{H^{2k}}^3+C\|\zeta\|_{H^{2k}}^3.
$$

Comparing with the ordinary differential equation
\begin{equation}\left\{\begin{split}&y'(t)=Cy^\frac{3}{2},\\
&y(0)=\|\omega_0\|_{H^{2k}}^2+\|\zeta_0\|_{H^{2k}}^2,
\end{split}\right.\end{equation}
where $\omega_0=\nabla\times\bf{u}_0$ and
$\zeta_0=\nabla\times\bf{H}_0$, then denoting by $T^*$ the blow up
time, it follows that  a priori estimates hold, for $T<T^*$,
\begin{equation}\begin{split}&\|\omega\|_{L^\infty(0,T;H^{2k})}+\|\zeta\|_{L^\infty(0,T;H^{2k})}\leq
C.
\end{split}\end{equation}

Thus, we have the following result.
\begin{theorem}
Let $\bf{u}_0\in V^{2k-1}\cap H^{2k+1}$ and $\bf{H}_0\in
V^{2k-1}\cap H^{2k+1}$, $k\geq 1$. Then there exist $T$ and
$C(\|\bf{u}_{0}\|_{H^{2k+1}},T)$, s.t.
\begin{equation}\begin{split}&\|\bf{u}\|_{L^\infty(0,T;H^{2k+1})}+\|\bf{H}\|_{L^\infty(0,T;H^{2k+1})}
\leq C.
\end{split}
\end{equation}\label{th1}
\end{theorem}

Taking the inner product $((\ref{om1}),
\partial_t\omega)_{H^{2k-1}}+((\ref{om2}),
\partial_t\zeta)_{H^{2k-1}}$, one obtains that
$\|\partial_t\omega\|_{L^2(0,T;H^{2k-1})}+\|\partial_t\zeta\|_{L^2(0,T;H^{2k-1})}\leq
C$. It follows that
$\|\partial_t\bf{u}\|_{L^2(0,T;H^{2k})}+\|\partial_t\bf{H}\|_{L^2(0,T;H^{2k})}\leq
C$.

According to  $\|\omega_0\|_{H^{2k}}\leq C$ and
$\|\zeta_0\|_{H^{2k}}\leq C$, then by equations (\ref{om1}) and
(\ref{om2}),
$\|\partial_t\omega|_{t=0}\|_{H^{2k-2}}+\|\partial_t\zeta|_{t=0}\|_{H^{2k-2}}\leq
C$

Similarly, applying operator $\partial_t\partial_{x,y,z}^{\alpha}$
to both sides of equations ($\ref{om1}$) and ($\ref{om2}$), for
$|\alpha|\leq 2k-2$,  and multiplying
$\partial_t\partial^\alpha_{x,y,z}\omega$ and
$\partial_t\partial^\alpha_{x,y,z}\zeta$, respectively,  we have,
\begin{equation}\begin{split}&\|\partial_t\bf{u}\|_{L^\infty(0,T;H^{2k-1})}+
\|\partial_t\bf{H}\|_{L^\infty(0,T;H^{2k-1})}\leq C.
\end{split}
\end{equation}

Thus, we can conclude
\begin{theorem}\label{th2}
Let the conditions of Theorem \ref{th1} be satisfied, then for
$s\leq k$
\begin{equation}\begin{split}&\|\partial^s_t\bf{u}\|_{L^\infty(0,T;H^{2k+1-2s})}
+\|\partial^s_t\bf{H}\|_{L^\infty(0,T;H^{2k+1-2s})}\leq C,\\&
\|\partial^{s+1}_t\bf{u}\|_{L^2(0,T;H^{2k-2s})}+
\|\partial^{s+1}_t\bf{H}\|_{L^2(0,T;H^{2k-2s})}\leq C.
\end{split}
\end{equation}
where $C=C(\|\bf{u}_0\|_{H^{2k+1}},T )$. \label{th2}
\end{theorem}

Then, the regularity of the solution of MHD equations (\ref{MS}) and
(\ref{bom}) is investigated,

\begin{theorem}
Let the conditions of Theorem \ref{th1} be satisfied. Then for
$s\leq k$, there exist a time $T$ depending on the initial date and
 unique classical solution of MHD equations (\ref{MS})
with boundary condition (\ref{bom}). In addition,
\begin{equation}\begin{split}&\|\partial^s_t\bf{u}\|_{L^\infty(0,T;H^{2k+1-2s})}
+\|\partial^s_t\bf{H}\|_{L^\infty(0,T;H^{2k+1-2s})}\leq C,
\\&
\|\partial^{s+1}_t\bf{u}\|_{L^2(0,T;H^{2k-2s})}+
\|\partial^{s+1}_t\bf{H}\|_{L^2(0,T;H^{2k-2s})}\leq C,
\end{split}
\end{equation}
where $C=C(\|\bf{u}_0\|_{H^{2k+1}},\|\bf{H}_0\|_{H^{2k+1}},T)$.
\label{th3}
\end{theorem}

\section{The vanishing viscosity limit}\label{section4}
This section focuses on the vanishing viscosity limit of the MHD
system for the flat boundary case.

\begin{theorem}\label{th4}Let the conditions of Theorem \ref{th1} be satisfied for $k\ge1$. Then as
$(\nu,\mu)\rightarrow0$, $(\bf{u},\bf{H})$ converge to the unique
solution $(\bf{u}^0,\bf{H}^0)$ of the ideal MHD system with the same
initial date in the sense
\begin{equation}\begin{split}(\bf{u},\bf{H})\rightarrow
(\bf{u}^0,\bf{H}^0) \text{ in }C(0,T;H^{2k}).
\end{split}
\end{equation}

\end{theorem}
\textbf{Proof.} It follows from Theorem \ref{th3} that
\begin{equation}\bf{u}(\nu,\mu),\bf{H}(\nu,\mu)\text{ is uniformly bounded in }
L^\infty(0,T;H^{2k+1}),\end{equation} and
\begin{equation}
\partial_t\bf{u}(\nu,\mu), \partial_t\bf{H}(\nu,\mu)\text{ is uniformly bounded in }
L^2(0,T;H^{2k}),
\end{equation}
for all $\nu>0$ and $\mu>0$. By the standard compactness result,
there exist a subsequence $\nu_k$ of $\nu$,  $\mu_k$ of $\mu$ and
vector functions $\bf{u}^0$ and $\bf{H}^0$, such that
\begin{equation}\begin{split}(\bf{u}(\nu_k,\mu_k),\bf{H}(\nu_k,\mu_k))\rightarrow
(\bf{u}^0,\bf{H}^0) \text{ in }C(0,T;H^{2k}),
\end{split}
\end{equation}
 as $(\nu,\mu)\rightarrow 0$. Passing to the limit, we can find
$(\bf{u}^0,\bf{H}^0)$ solves the ideal MHD equations (\ref{IM}) and
(\ref{boi}). Together with the uniqueness of the strong solution of
the ideal MHD systems, we then show the convergence of whole
sequence. \Box

Now, we present the convergence rate.
\begin{theorem}Let the conditions of Theorem \ref{th1} be satisfied
for $k\geq 1$. Then,
\begin{equation}\|\bf{u}-\bf{u}^0\|_{L^\infty(0,T;H^{2k-1})}
+\|\bf{H}-\bf{H}^0\|_{L^\infty(0,T;H^{2k-1})}\leq C\nu+C\mu.
\end{equation}\label{th5}
\end{theorem}
\textbf{Proof.} Set $\omega^0=\nabla\times\bf{u}^0$ and
$\zeta^0=\nabla\times\bf{H}^0$. Recalling Lemma \ref{lem2} and
Theorem \ref{th4}, one obtains
\begin{equation}\begin{split}&\bf{u}^0(t)\in X^{2k-1}\cap H^{2k+1}, \bf{H}^0(t)\in X^{2k-1}\cap
H^{2k+1},\\
&\omega^0(t)\in X^{2k-2}\cap H^{2k}, \zeta^0(t)\in X^{2k-2}\cap
H^{2k}.\end{split}
\end{equation}

Set $\widetilde{\bf{u}}=\bf{u}-\bf{u}^0$,
$\widetilde{\omega}=\nabla\times\widetilde{\bf{u}}$,
$\widetilde{\bf{H}}=\bf{H}-\bf{H}^0$,
$\widetilde{\zeta}=\nabla\times\widetilde{\bf{H}}$. We can find
$\widetilde{\bf{u}}, \widetilde{\omega}$ solve
\begin{equation}\begin{split}&\quad\partial_t\widetilde{\omega}-\nu\triangle
\widetilde{\omega}+(\bf{u}\cdot\nabla)\widetilde{\omega}
+(\widetilde{\bf{u}}\cdot\nabla)\omega^0-(\bf{H}\cdot\nabla)\widetilde{\zeta}
-(\widetilde{\bf{H}}\cdot\nabla)\zeta^0\\
&+F(D\bf{u},D\widetilde{\bf{u}})+F(D\widetilde{\bf{u}},D\bf{u}^0)
-F(D\bf{H},D\widetilde{\bf{H}})-F(D\widetilde{\bf{H}},D\bf{H}^0)=
\nu\triangle\omega^0,
 \end{split}\label{om3}
\end{equation}
\begin{equation}\begin{split}&\quad\partial_t\widetilde{\zeta}-\mu\triangle
\widetilde{\zeta}+(\bf{u}\cdot\nabla)\widetilde{\zeta}
+(\widetilde{\bf{u}}\cdot\nabla)\zeta^0-(\bf{H}\cdot\nabla)\widetilde{\omega}
-(\widetilde{\bf{H}}\cdot\nabla)\zeta^0\\
&+F(D\bf{u},D\widetilde{\bf{H}})+F(D\widetilde{\bf{u}},D\bf{H}^0)
-F(D\bf{H},D\widetilde{\bf{u}})-F(D\widetilde{\bf{H}},D\bf{u}^0)=
\mu\triangle\omega^0,
 \end{split}\label{om4}
\end{equation}
Noting that  $\widetilde{\omega}^0\in V^{2k-2}$,
$\widetilde{\zeta}^0\in V^{2k-2}$, $\|\triangle
\omega\|_{H^{2k-2}}\leq C$ and $\|\triangle \zeta\|_{H^{2k-2}}\leq
C$, the same argument in proof of Theorem \ref{th1} can be followed.
Taking the inner products
$((\ref{om3}),\widetilde{\omega})_{H^{2k-2}}+((\ref{om4}),\widetilde{\zeta})_{H^{2k-2}}$,
and the desired result can be obtained. $\Box$

There is a gap between $(\partial^{2k}_n\omega^0_\tau,
\partial^{2k}_n\zeta^0_\tau)$ and 0. In  other words, $\omega^0$
and $\zeta^0$ are not in $V^{2k}$.  Assuming
$\|\partial^{2k}_n\omega_\tau^0\|_{L^\infty(0,T;C^2(\partial\Omega))}\leq
C$,
$\|\partial^{2k}_n\zeta_\tau^0\|_{L^\infty(0,T;C^2(\partial\Omega))}\leq
C$,
$\|\partial_t\partial^{2k}_n\omega_\tau^0\|_{L^2(0,T;C^2(\partial\Omega))}\leq
C$,
$\|\partial_t\partial^{2k}_n\zeta_\tau^0\|_{L^2(0,T;C^2(\partial\Omega))}\leq
C$,  by Lemma \ref{lem3},  there exist  $\bf{v}^{i}\in V^{2k-1}$,
$\chi^{i}=\nabla\times\bf{v}^\nu\in V^{2k-2}$ ($i=1,2$), s.t.
$\partial^{2k}_n\chi^1_\tau=-\partial^{2k}_n\omega^0_\tau$,
$\partial^{2k}_n\chi^2_\tau=-\partial^{2k}_n\zeta^0_\tau$ on
$\partial\Omega$,
$\|z^j(1-z)^j\partial^{2k+1}_n\chi^i\|_{L^\infty(0,T;L^{2})}\leq
C\nu^{\frac{2j-1}{4}}$ , and further
$\|\partial_t\chi^i\|_{L^2(0,T;H^{2k})}\leq C\nu^{\frac{1}{4}}$,
$\|\partial_\tau\chi^i\|_{L^\infty(0,T;H^{2k})}\leq
C\nu^{\frac{1}{4}}$, $i=1,2$, $j=0,1$.

Set $\widehat{\bf{u}}=\bf{u}-\bf{u}^0-\bf{v}^1$,
$\widehat{\bf{H}}=\bf{H}-\bf{H}^0-\bf{v}^2$
$\widehat{\omega}=\nabla\times\widehat{\bf{u}}=\omega-\omega^0-\chi^1$,
$\widehat{\zeta}=\nabla\times\widehat{\bf{H}}=\zeta-\zeta^0-\chi^2$.

From
equations (\ref{om3}) and (\ref{om4}), one obtains,
\begin{equation}\begin{split}&\quad\partial_t\widehat{\omega}-\nu\triangle
\widehat{\omega}+(\bf{u}\cdot\nabla)\widehat{\omega}
+(\widehat{\bf{u}}\cdot\nabla)\omega^0-(\bf{H}\cdot\nabla)\widehat{\zeta}
-(\widehat{\bf{H}}\cdot\nabla)\zeta^0\\
&+F(D\bf{u},D\widehat{\bf{u}})+F(D\widehat{\bf{u}},D\bf{u}^0)
-F(D\bf{H},D\widehat{\bf{H}})+F(D\widehat{\bf{H}},D\bf{H}^0)\\
&= \nu\triangle\omega^0-
\partial_t\chi^1+\nu\triangle\chi^0-(\bf{u}\cdot\nabla)\chi^1
-(\bf{v}^1\cdot\nabla)\omega^0+(\bf{H}\cdot\nabla)\chi^2
+(\bf{v}^2\cdot\nabla)\zeta^0\\
&-F(D\bf{u},D\bf{v}^1)-F(D\bf{v}^1,D\bf{u}^0)
+F(D\bf{H},D\bf{v}^2)+F(D\bf{v}^2,D\bf{H}^0) ,
 \end{split}\label{om5}
\end{equation}
\begin{equation}\begin{split}&\quad\partial_t\widehat{\zeta}-\mu\triangle
\widehat{\zeta}+(\bf{u}\cdot\nabla)\widehat{\zeta}
+(\widehat{\bf{u}}\cdot\nabla)\zeta^0-(\bf{H}\cdot\nabla)\widehat{\omega}
-(\widehat{\bf{H}}\cdot\nabla)\zeta^0\\
&+F(D\bf{u},D\widehat{\bf{H}})+F(D\widehat{\bf{u}},D\bf{H}^0)
-F(D\bf{H},D\widehat{\bf{u}})-F(D\widehat{\bf{H}},D\bf{u}^0)\\&=
\mu\triangle\omega^0-\partial_t\chi^2+\mu\triangle
\chi^2-(\bf{u}\cdot\nabla)\chi^2
-(\bf{v}^1\cdot\nabla)\zeta^0+(\bf{H}\cdot\nabla)\chi^1
+(\bf{v}^2\cdot\nabla)\zeta^0\\
&-F(D\bf{u},D\bf{v}^2)-F(D\bf{v}^1,D\bf{H}^0)
+F(D\bf{H},D\bf{v}^1)+F(D\bf{v}^2,D\bf{u}^0),
 \end{split}\label{om6}
\end{equation}

Then, taking the inner products
$((\ref{om5}),\widehat{\omega})_{H^{2k}}+((\ref{om6}),\widehat{\zeta})_{H^{2k}}$,

Note that
$$\partial^{2k-1}_n\widehat{\omega}_n=0,
\partial^{2k}_n\widehat{\omega}_\tau=0, \partial^{2k-1}_n\widehat{\zeta}_n=0,
\partial^{2k}_n\widehat{\zeta}_\tau=0\text{ on }\partial\Omega,
$$
$$|\nu\int\triangle \partial_{n}^{2k}\chi^1\cdot \partial_{n}^{2k}\widehat{\omega}^\nu|=
|\nu\int\nabla \partial_{n}^{2k}\chi^1\cdot
\nabla\partial_{n}^{2k}\widehat{\omega}^\nu|\leq
\frac{\nu}{8}\|\nabla\partial_{n}^{2k}\widehat{\omega}^\nu\|_{L^2}^2+C\nu^{\frac{1}{2}},$$
$$|\nu\int\triangle \partial_{n}^{2k}\omega^0\cdot \partial_{n}^{2k}\widehat{\omega}^\nu|=
|\nu\int\nabla \partial_{n}^{2k}\omega^0\cdot
\nabla\partial_{n}^{2k}\widehat{\omega}^\nu|\leq
\frac{\nu}{8}\|\nabla\partial_{n}^{2k}\widehat{\omega}^\nu\|_{L^2}^2+C\nu,$$
$$\|(\bf{u}\cdot\nabla)\partial^{2k}_n\chi^1\|_{L^2}\leq
\|(\bf{u}_\tau\cdot\nabla_\tau)\partial^{2k}_n\chi^1\|_{L^2}+
\|\frac{u_3}{z(1-z)}\|_{L^\infty}\|z(1-z)\partial_n^{2k+1}\chi^1\|_{L^2}\leq
C\nu^{\frac{1}{4}}.$$

And it follows in the same manner, that
$$\|\widehat{\omega}\|_{L^\infty(0,T;H^{2k})}+\|\widehat{\zeta}\|_{L^\infty(0,T;H^{2k})}\leq C\nu^{\frac{1}{4}}+C\mu^{\frac{1}{4}}.$$

The following result is concluded,
\begin{theorem}Let the conditions of Theorem \ref{th1} be satisfied for $k\geq 1$.
Assume the solution $(\bf{u}^0,\bf{H}^0)$ of the ideal MHD equations
(\ref{IM}) and (\ref{boi}) satisfy
$\|\partial^{2k}_n\omega_\tau^0\|_{L^\infty(0,T;C^2(\partial\Omega))}\leq
C$,
$\|\partial^{2k}_n\zeta_\tau^0\|_{L^\infty(0,T;C^2(\partial\Omega))}\leq
C$,
$\|\partial_t\partial^{2k}_n\omega_\tau^0\|_{L^2(0,T;C^2(\partial\Omega))}\leq
C$,
$\|\partial_t\partial^{2k}_n\zeta_\tau^0\|_{L^2(0,T;C^2(\partial\Omega))}\leq
C$, then,
\begin{equation}\|\bf{u}-\bf{u}^0\|_{L^\infty(0,T;H^{2k+1})}+\|\bf{H}-\bf{H}^0\|_{L^\infty(0,T;H^{2k+1})}\leq
C\nu^{\frac{1}{4}}+C\mu^{\frac{1}{4}}.
\end{equation}\label{th6}
\end{theorem}

Finally, we give two remarks.

\begin{remark} If the conditions of Theorem \ref{th5} are all
satisfied, then for $s\leq k-1$,
\begin{equation}\|\partial_t^s\bf{u}-\partial^t_s\bf{u}^0\|_{L^\infty(0,T;H^{2k-1-2s})}
+\|\partial_t^s\bf{H}-\partial^t_s\bf{H}^0\|_{L^\infty(0,T;H^{2k-1-2s})}\leq
C\nu+C\mu.
\end{equation}
\end{remark}
\begin{remark}
If the conditions of Theorem \ref{th6} are all satisfied, then for
$s\leq k$,
\begin{equation}\|\partial^s_t\bf{u}-\partial_t^s\bf{u}^0\|_{L^\infty(0,T;H^{2k+1-2s})}+
\|\partial^s_t\bf{H}-\partial_t^s\bf{H}^0\|_{L^\infty(0,T;H^{2k+1-2s})}
\leq C\nu^{\frac{1}{4}}+ C\mu^{\frac{1}{4}}.
\end{equation}
\end{remark}


\bibliographystyle{amsplain}

\end{document}